\documentclass{article}

\newtheorem{thm}{Theorem}[section]
\newtheorem{lem}[thm]{Lemma}
\newtheorem{dfn}[thm]{Definition}

  \ifx\pdfoutput\undefined
   \usepackage[dvips]{graphicx}
   \else
   \usepackage[pdftex]{graphicx}
   \fi
   \usepackage{epstopdf}
 \DeclareGraphicsExtensions{.pdf,.gif,.jpg}
 
 \usepackage{amsfonts}
  \usepackage{amsmath}
 
\begin{document}

\title{A pedagogical history of compactness}
\author{Manya Raman-Sundstr\"{o}m \thanks{This paper is based upon my masters thesis~\cite{rama} at UC Berkeley. The paper has benefitted, at different stages, from the generous help of the following people:  Hendrik Lenstra, Hans Wallin, Lucien Le Cam, Johan de Jong, Jeremy Gray, Klas Markstr\"{o}m, Victor Falgas-Ravry, Edouard Servan-Schreiber, and Lars-Daniel \"{O}hman. I'm also grateful to three anonymous reviewers and several librarians, including Mikael R\r{a}gstedt at the Mittag-Leffler Institute, who helped me track down original sources. This project was funded in part by a stipend from Wenner-Gren Foundation. Of course all mistakes are my own.}}

\maketitle


\begin{flushright}
\parbox{2.8 in}{Modern mathematics tends to obliterate history:  
each new school rewrites the foundations of its subject 
in its own language, which makes for fine logic but poor pedagogy.}

\emph{R. Hartshorne}
\end{flushright}


\section{Why study the history of compactness?}

Compactness has come to be one of the most important and useful notions in advanced
mathematics. It can also be seen as a kind of a gate-keeper topic: if an undergraduate mathematics student does not understand compactness, whatever it really means to understand, it is unlikely that he or she will be able to do higher level mathematics. However, for whatever reasons, when we teach compactness (and other topics of similar importance) we tend to do so without motivation, leaving students on their own to figure out, if they ever do, how various definitions and theorems relate to each other and why they take the specific forms they do.

This paper is an attempt to fill in some of the information that the standard textbook treatment of compactness leaves out. It is not a historical article, per se, but a synthesis of historical documents with an eye towards clarifying the main ideas related to compactness. In particular, the paper discusses the origins and development of both open-cover and sequential compactness, how and why open-cover compactness came to be favored, and some modern developments involving nets and filters.

A list of terms related to compactness is given in the Appendix. Since the terms have changed names at various points in history, the list can be useful for keeping the concepts straight. In the main text we will use a combination of historical and modern formulations of the main definitions, lemmas, and theorems (choosing the ones best for readability, with originals in footnotes.)\\

\section{Possible motivations for compactness}
Compactness grew out of one of the most productive periods of mathematical activity. 
In middle to late nineteenth century Europe, advanced mathematics began to take the
form in which we know it today. In the background was Cantor's work establishing the
beginning of a systematic study of set theory and point-set topology (though Cantor himself
turned his interests to transfinite sets. The significance of Cantor's work for topology should
be credited to Poincar\'{e}\footnote{Thanks to Jeremy Gray for this comment.}). Also, many
mathematicians---including Weierstrass, Hausdorff and Dedekind---were worried about the
foundations of mathematics and began to make rigorous many of the ideas that had for
centuries been taken for granted. While some of the nineteenth century work can be
traced to mathematical concerns of the early Greeks, the level of rigor and
abstraction reflects a revolution in mathematical thought.

It is in this context that we will discuss some specific problems that appear to have
motivated the concept of compactness. In particular, we will discuss the influence of
the study of properties of closed, bounded intervals of real numbers (which I will
denote $[a,b]$), spaces of continuous functions, and solutions to differential
equations.

\subsection{Properties of $[a,b]$}

In mid to late nineteenth century, mathematicians began to really understand and specify essential properties of the real line. There were essentially two characterizations that were developed during this time. One characterization, developed by Bolzano and Weierstrass among others, grew out of the study of functions defined on sequences of real numbers. The other characterization, which grew out of work by Heine, Borel, and Lebesgue, was based on topological features, such as the covering of sets by open neighborhoods. We will examine both of these characterizations in more detail.

The origin of sequential compactness is often traced to a theorem, proved rigorously by Weierstrass in 1877, which concerns the behavior of continuous functions defined on closed, bounded intervals of the real line.\footnote{This date refers to one of the earliest publications of the theorem, see~\cite{weie}. However it is likely that Weierstrass actually proved it years before and disseminated it orally, via lectures, perhaps ten years earlier.} The following statement of the theorem comes from Fr\'{e}chet, who referred to this theorem as a result of Weierstrass.

\begin{thm}[Weierstrass] Each function continuous in a limited [equivalent to modern-day ``closed and bounded''] interval attains there at least once its maximum.\footnote{From the original French: Weierstrass a en effete d\'{e}montr\'{e} que toute fonction continu\'{e} dans un intervalle limit\'{e} y atteint au moins une fois son maximum.~\cite[p. 848]{fre1}}
\label{thm:WMT}
\end{thm}

\begin{figure}[h]
\label{fig1}
 � \begin{center}
 � � \includegraphics[width=2in]{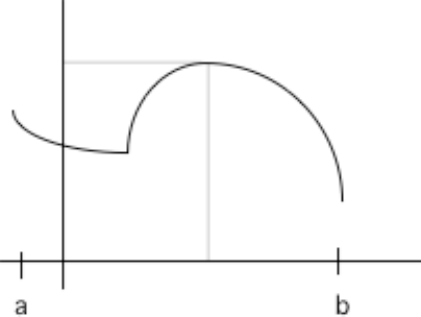}
 � \end{center}
 � \caption{A continuous function on $[a,b]$.}
  \end{figure}

Fr\'{e}chet, who defined sequential compactness in his 1906 thesis, said that his definition came from his desire to generalize this theorem to abstract topological spaces~\cite[p. 244]{tay1}. Weierstrass's theorem owes its essential ideas to Bolzano in 1817, who working in relative isolation, both politically and mathematically, in Bohemia, remarkably stated and proved the following:  

\begin{lem}[Bolzano]If a property $M$ does not apply to all values of a
 variable quantity $x$, but to all those that are {\em smaller} than a certain $u$, there is
always a quantity $U$ which is the greatest of those of which it can be asserted that
all smaller $x$ possess the property $M$.\footnote{From the original German: Wenn eine Eigenschaft $M$ nicht allen Werthen einer ver\"{a}nderlichen Gr\"{o}sse $x$, wohl aber allen, die kleiner sind, als ein
gewisser $u$, zuk\"{o}mmt: so gibt es allemahl eine Gr\"{o}sse $U$, welche die
gr\"{o}sste derjenigen ist, von denen behauptet werden kann, dass alle
kleineren $x$ die Eigenschaft $M$ besitzen.~\cite[p. 41]{bolz}}
\label{lem:Bolzano}
\end{lem}

This lemma, today called the greatest lower bound property for real numbers, was somewhat of a breakthrough in the conceptualization of real numbers. The proof of this lemma provided the first real account of the limiting process, and was used to prove what we now call the Intermediate Value Theorem: if $f$ is continuous on $[a,b]$ with $f(a)<0$ and $f(b)>0$, then for some $x$ between $a$ and
$b$, $f(x)$ will be exactly $0$.  

The idea behind Bolzano's proof of the lemma was to use interval bisection, that is to narrow in on the least upper bound by throwing away points of the set that were below it. This iterative process was essentially 
the same process used in Weierstrass's proof of the maximum value
theorem (see ~\cite[p.953]{klin}). In particular, Bolzano's lemma allowed Weierstrass
to prove that every bounded infinite set of real numbers has a limit point. It is
this property that Fr\'{e}chet used when he generalized Weierstrass's theorem to
abstract spaces. We now know this property as the {\em Bolzano-Weierstrass property},
or {\em limit-point compactness}.

While Bolzano and Weierstrass were trying to characterize properties of the real line
in terms of sequences, other mathematicians, such as Borel and Lebesgue, were trying to
characterize it in terms of open covers. Borel proved the following lemma in his 1894
thesis:  

\begin{lem}[Borel] If on a line one has an infinite number of subintervals, such that every point of the line is interior to at least one of the intervals then one can determine effectively a bounded number of intervals from among the given intervals that have the same property (every point of the line is interior to at least one of them.)\footnote{From the original French:  Si l'on a sur une droite une infinit\`{e} d'intervalles partiels, tels que tout point de la droite soit int\'{e}rieur \`{a} l'un au moins des intervalles, on peut d\'{e}terminer effectivement un nombre limit\'{e} d'intervalles choisis parmi les intervalles donn\'{e}s et ayant la m\^{e}me propri\'{e}t\'{e} (tout point de la droite est int\'{e}rieur \`{a} au moins l'un d'eux).~\cite[p. 51]{bore}}

\end{lem}

Here a line means a bounded interval. It turns out that Borel's approach was similar to the approach Heine used in 1872 to prove that a continuous function on a closed interval was uniformly continuous~\cite[p. 188]{hein}. This
theorem was first proven by Dirichlet in his lectures of 1852, with a more explicit
use of coverings and subcoverings than in Heine's theorem~\cite[p. 91]{duga}. However
Dirichlet's notes were not published until 1904, which might explain why he does not get
credit for the generalized version of the Borel lemma (now referred to as Borel Theorem). The reason that Heine's name
is attached to the theorem is that Sch\"{o}nflies, a student of Weierstrass, noticed
the connection between Heine's work and Borel's~\cite[p. 51]{schoe}. The generalized theorem, which is
now commonly called the Heine-Borel theorem,\footnote{For an accessible history of this theorem, along with discussion of a number of different original statements of this theorem, see~\cite{andr}.} with modern language and notation, is:

\begin{thm}[Heine-Borel]  A subset of $\mathbb{R}$ is compact iff it is closed and bounded.
\end{thm}

While Heine is credited with a theorem he did not prove, it appears that Cousin was
largely overlooked for a lemma he did prove. In 1895, he generalized the Borel
lemma to arbitrary covers. The following is often referred to as Cousin's Theorem,
but it appears in the original as a lemma. The plane YOX below is just $\mathbb{R}^2$, and the region $S$ would, in today's language,
be described as closed and bounded.

\begin{lem}
[Cousin] In the plane YOX let S be a connected area bounded by a closed contour, simple or complex; one supposes that at each point of S or its perimeter there is a circle, of non-zero radius, having this point as its centre; it is then always possible to subdivide S into regions, finite in number and sufficiently small for each one of them to be entirely inside a circle corresponding to a suitably chosen point in S or on its perimeter.\footnote{From the original French: Soit, sur le plan YOX, une aire connexe $S$ limit\'{e}e par un contour ferm\'{e} simple ou complexe; on suppose qu'\`{a} chaque point de $S$ ou de son p\'{e}rim\`{e}tre correspond un cercle, de rayon non nul, ayant ce point pour centre: il est alors toujours possible de subdivider $S$ en r\'{e}gions, en nombre fini et assez petites pour que chacune d'elles soit compl\`{e}tement int\'{e}rieure au cercle correspondant \`{a} un point convenablement choisi dans $S$ ou sur son p\'{e}rim\`{e}tre.~\cite[p. 22]{cous}}
\end{lem}

In other words, if for every point of a closed, bounded region, there corresponds a finite neighborhood, then the region can be divided into a finite number of subregions such that each subregion is contained in a circle having its center in the subregion.\footnote{Note that the original definition was formulated without the term 'neighborhood,' or 'voisinage' in French. The idea of neighborhood was obviously around during Cousin's time, but was not used consistently. Formal definitions of the term can be found in~\cite{fre2} and~\cite{haus}).} Cousin's theorem is generally attributed to Lebesgue, who was said to be aware of the result in 1898 and published his proof in 1904~\cite{lebe}\footnote{Cited in ~\cite[p. 29]{hild}.}. The Lebesgue lemma is considered itself to be an important consequence of compactness.

While there is some debate over who was really responsible for the ideas and proofs,
the idea that any closed, bounded subset of $\mathbb{R}$ has the open-cover
property (sometimes called the Borel-Lebesgue property) was known when Fr\'{e}chet first defined compactness formally. 

\subsection{Spaces of continuous functions}

A second motivation for the notion of compactness was the study of abstract topological spaces such as spaces of continuous functions,\footnote{We could just as well take $C^0$ on any set $\Omega$.} $C^0[a,b]$. In $C^0[a,b]$, points are functions (whereas in $[a,b]$ points are real numbers). The properties of $[a,b]$ alone might not have been seen as important to generalize if it weren't the case that these properties seemed to be important in more abstract spaces as well. However, it turned out that infinite dimensional spaces like $C^0[a,b]$ were not as well-behaved as finite dimensional spaces like $\mathbb{R}$. For instance, closed, bounded
subsets of continuous functions on $\mathbb{R}$ do not necessarily have the Bolzano-Weierstrass or open-cover property. The work in this area was done by Ascoli and Arzel\`{a} in the last decades of the 1800's.

The following example illustrates that a closed, bounded subset
of continuous functions on $\mathbb{R}$ is not, in our modern language, sequentially compact. 
Consider $B$, the set of continuous functions, $f$, defined on $[0,1]$ with $\|f\| \leq1$.  
(This is the closed unit ball in $C^0[a,b]$ and $\|$ $\|$ is the sup
norm.) We will show that there is a sequence in $B$ that does not have a convergent
subsequence. Let $f_n(x)=x^n$. This sequence lies in $B$, but we cannot find a subsequence
that converges to a function in $C^0[a,b]$. Suppose to the contrary $f$ is such a function. Then

\[f(x) =\lim_{k \rightarrow \infty} f_{n_k}(x) \] 

which would imply that 

\[f(x)=\left\{ \begin{array}{ll}
               0  & \mbox{if $x<1$}\\
		               1  & \mbox{if $x=1$}
               \end{array}.
       \right. \]

Since $f$ is a discontinuous function, it is not in $C^0[a,b]$. Hence the sequence
$f_n(x)$ has no convergent subsequence.

The problem in this example comes from how the functions converge. If convergence
means pointwise convergence, we do not get behavior analogous to that of, say,
sequences in closed unit balls of $\mathbb{R}$. In order to avoid this problem, Ascoli
introduced the notion of equicontinuity~\cite[p. 566]{asc1}.\footnote{See also~\cite[p. 27]{bour}.} Equicontinuity
requires functions to converge to a limit all at once instead of pointwise. A set $E$
is {\em equicontinuous} iff for all $\epsilon > 0$ there exists a $\delta>0$  such
that $|s-t| < \delta$ and $f \in E$ imply $|f(s)-f(t)|< \epsilon$.

The Arzel\`{a}-Ascoli theorem, in modern language, then states: 

\begin{thm} [Arzel\`{a}-Ascoli]
 Any bounded equicontinuous sequence of functions in
$C^0[a,b]$ has a uniformly convergent subsequence.\footnote{From the original Italian: La condizione necessaria e sufficiente affinch\'{e} una successione data di funzioni $u_1(x), u_2(x), ... u_n(x),...$ abbia una funzione limite, nel senso detto sopra, \'{e} che, preso un numero positivo $\sigma$ piccolo a piacere, si possa sempre determinare un numero intero corrispondente $m_\sigma$ tale che per ogni $x$, nell'intervallo $a ... b$, si abbia qualunque sia $p$ intero positivo, $\mid u_{m_\sigma(x)} - u_{m_\sigma+p} (x) \mid < \sigma$. Translated to English: The necessary and sufficient conditions that a given sequence of functions $u_1(x), u_2(x), ... u_n(x),...$ defined on an interval $a ... b$ may converge to a limiting function is that, given an arbitrarily small positive number $\sigma$ there can always be determined a corresponding integer $m_\sigma$ such that for all values of $x$ in the interval $a ...b$, and for all positive integers $p$, $\mid u_{m_\sigma(x)} - u_{m_\sigma+p} (x) \mid < \sigma$~\cite[p. 226]{asc2}.}
\end{thm}

Using modern terminology we can state a consequence of this theorem, analogous to the Heine-Borel theorem:

\begin{thm}
A subset of $C^0[a,b]$ is compact iff it is closed, bounded and equicontinuous.
\end{thm}

Ascoli proved the sufficiency of this condition in 1884~\cite[p.567]{asc1} and Arzel\`{a} the necessity in 1889~\cite[p. 345]{arze} (with a clearer proof presented in 1894~\cite[p. 226]{asc2}). This generalization of Bolzano-Weierstrass's theorem (although not stated in terms of compactness) was apparently well known after 1880. Moreover, Hilbert seems to have discovered this
``compactness" property independently and published it in 1900~\cite[p. 82]{die1}. It is 
unclear whether Arzel\`{a} and Ascoli themselves were aware of how their work
was connected with compactness, but Fr\'{e}chet's work was influenced by theirs~\cite[p. 255]{tay1}.


\subsection{Solutions to differential equations}

A third motivation for the notion of compactness came from the desire to find
solutions to differential equations. Peano, a contemporary of Arzel\`{a} and Ascoli as
well as a fellow Italian, realized that the Arzel\`{a}-Ascoli theorem might be useful
for demonstrating the existence of such solutions. He searched for solutions by
making a sequence of approximations. He then used what we now call compactness to
show that there will be a subsequence that converges to a limit (which will be the
solution to the differential equation).

To this end, Peano proved the following theorem in 1890:

 \begin{thm} [Peano]

Suppose we are given a system of differential equations in normal form:
\begin{align*}
 \frac{dx_1}{dt}&= \varphi_1(t, x_1, \ldots ,x_n),\\
\ldots & \ldots \\
\frac{dx_n}{dt}&= \varphi_n(t, x_1, \ldots ,x_n)
\end{align*}
where the functions $\varphi_1, \ldots, \varphi_n$ are continuous in a neighbourhood of $(b, a_1, \ldots, a_n)$. [Then there exists] an interval $(b,b')$, and $n$ functions of $t$ from this interval, $x_1, \ldots, x_n$, satisfying our system of equations and evaluating to $a_1, \ldots, a_n$ respectively at $t=b$.\footnote{From the original French: Soit donn\'{e} le syst\`{e}me d'\'{e}quations diff\'{e}rentielles, ramen\'{e} \`{a} la forme normale:$ \frac{dx_1}{dt}= \varphi_1(t, x_1, ... ,x_n), ....  \frac{dx_n}{dt}= \varphi_n(t, x_1, ... ,x_n)$, o\`{u} les $ \varphi_1, ..., \varphi_n$ sont des fonctions continues aux environs de $t = b, x_1 = a_1,..., x_n = a_n.$ [Alors il existe] un intervalle $(b,b')$, et, dans cet intervalle, $n$ fonctions $x1  \ldots x_n$ de $t$, qui satisfont  aux \'equations donn\'ees, et qui, pour $t=b$, prennent les values $a_1 \ldots a_n$.~\cite[p.182]{pean}}

\end{thm}

While it is not clear if Fr\'{e}chet was aware of this application, it does seem to be the case that applications for the notion of
compactness were known before the term was formally defined.  



\section{Developing the definition}

We will trace below the development of the two central notions of compactness
discussed above, those stemming from sequences and open covers of real numbers. 
Again, it is useful to know something about the climate of the mathematics community
at the time of these historical developments. We will focus on the contributions of
only a few central people, but there was actually a large community of mathematicians
who were developing ideas that are now the foundations for analysis and topology. 
Many of these mathematicians were in very close contact with each other, so it is difficult 
to tease apart their contributions. Among them, in France, were Hadamard, Lebesgue, and Fr\'{e}chet; in Russia, Alexandroff\footnote{I will use the spelling Alexandroff, rather than the sometimes used Alexandrov, since that 
was the spelling he preferred.} and Urysohn; in Germany, Hausdorff, Hilbert, Sch\"{o}nflies, and Cantor; in Hungary, F. Riesz; in the Netherlands, Brouwer; and in the U.S., Chittenden, Hedrick, and Moore.  

We will start with the work of Fr\'{e}chet, who coined the term ``compact'' and gave
definitions for what we now know as countable and sequential compactness. We will
then briefly discuss contributions by Alexandroff and Urysohn who developed and stated
what we now call open-cover compactness, or simply compactness. We will show why open-cover and sequential compactness are
not equivalent in abstract topological spaces, providing motivation for a formulation of compactness 
in terms of nets and filters which is analogous to sequential compactness.

\subsection{Fr\'{e}chet: Countable and limit-point compactness}

While Fr\'{e}chet was influenced by many contemporaries and predecessors, it seems he
deserves credit as the father of compactness. It was Fr\'{e}chet who gave the concept a name, 
in a paper \cite{fre1} leading to his 1906 doctoral thesis. Fr\'{e}chet also defined metric spaces for the first time, though not using that term,
 and made inroads into functional analysis, thus providing a context for which the importance of 
 compactness became clear.

Fr\'{e}chet was a mathematician of big ideas. He preferred definitions that had an
intuitive feel rather than analytic power. This preference can be seen in \cite[p. 849]{fre1} in which he defined a notion of compactness, introducing first a definition of what we now call countable compactness, using nested intersections, before introducing a characterisation
using limit points. 

In Fr\'{e}chet's thesis, he considered three kinds of spaces, which he called L-class,
V-class, and E-class. L-classes were the most general, in which a notion of sequential
compactness was defined.  E-classes, which we now call metric spaces, and V-classes,\footnote{The letter V comes from the French {\em voisinage} meaning neighborhood. L stands for  {\em limite}, or limit.  E stands for {\em \'{e}cart}, or gap, as in the non-zero distance between two points in a metric space.} a metric space with a weak version of the triangle inequality, were less general, but easier to work with. The goal was to define compactness for L-classes, but this turned out unsuccessful because sequential compactness did not have all the properties needed to generalize to abstract topological spaces (more on this in section 3.2). Fr\'{e}chet focused instead on the V- and E-classes, in which notions of modern-day compactness and sequential or limit-point compactness were equivalent. The following definition was given for E-classes.

\begin{dfn}
A set $E$ is called compact if, whenever ${E_n}$ is a sequence of nonempty, closed
subsets of $E$ such that $E_{n + 1}$ is a subset of ${E_n}$ for each $n$, there is at least
one element that belongs to all of the $E_n$'s.\footnote{From the original French:
Si on consid\`{e}re une suite d'ensembles $E_1, E_2, ..., E_n, ...$ form\'{e}s d'\'{e}l\'{e}ments d'un m\^{e}me ensemble compact $E$, chacun ferm\'{e}, contenu dans le pr\'{e}c\'{e}dent et poss\'{e}dent au moins un \'{e}l\'{e}ment, il y a n\'{e}cessairement un \'{e}l\'{e}ment commun \`{a} tous ces ensembles.~\cite[p.7]{fre2}}
\end{dfn}

The exact nature of Fr\'{e}chet's intuition for this definition is unclear, but there might be two features of compact sets that he wanted to capture.  The first is a sense of boundedness. The nested intersection property allows us to easily rule out sets that have tails running to $\infty$. For instance, we can show that $\mathbb{R}$ is not
compact. Let $E_n=[n,\infty)$. Each $E_n$ is closed because it contains the point $n$, and
clearly $E_{n+1} \subset E_n$. However, the infinite intersection of these intervals is
empty, so $\mathbb{R}$ is not compact. 

\begin{figure}[h]
\label{fig4}
 \begin{center}
\includegraphics[width=3in]{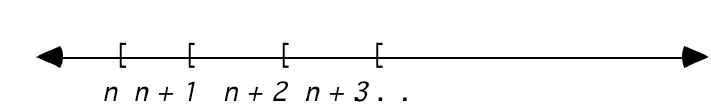}
 \end{center}
 \caption{Nested tails.}
 \end{figure}

A second feature the nested intersection definition allows us to quickly see is that
sets that have ``holes'' are not compact. For instance, we can see that $X = [a,b]$ is
compact and $Y = [a,b) \cup (b,c]$ is not.\footnote{Of course one could construct a similar 
example {\it with} a hole that is compact such as $Z = [a,b] \cup [c,d]$, but the example illustrates what could have been the intuition behind the definition.} In the latter case, consider $E_i=[a_i,b) \cup (b,c_i]$  where $a_{i+1}>a_i$ and $a_i \rightarrow b$, $c_{i+1}<c_i$, and $c_i \rightarrow b$ . These sets are clearly nested and are closed in $Y$, but the infinite intersection of those intervals is empty. Hence $Y$ is not compact. 

\begin{figure}[h]
\label{fig3}
 \begin{center}
 \includegraphics[width=3in]{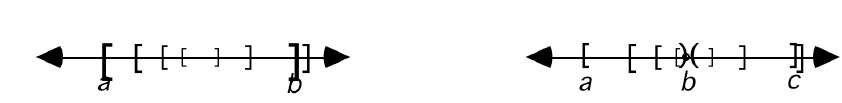}
 \end{center}
 \caption{Nested intersections.}
 \end{figure}

While Fr\'{e}chet preferred his intuitive definition involving nested intersections~\cite[p. 430]{pier}, he realized the importance of also providing a more useful, if less intuitive, definition. Below is another definition from Fr\'{e}chet, which uses the Bolzano-Weierstrass property. This definition applies to V- and E-classes where limit-point, countable, and sequential compactness are equivalent. Note that for Fr\'{e}chet, a compact set need not be closed.  So for his subsequent definitions and theorems, he often needed to require that a set be both compact and closed.

\begin{dfn}
We will say that a set is [relatively limit-point] compact if it contains only a finite number of points or if every one of its infinite subsets gives rise to at least one limit point.\footnote{From the original French: Nous dirons qu'un ensemble est {\it compact} lorsqu'il ne comprend qu'un nombre fini d'\'{e}l\'{e}ments ou lorsque toute infinit\'{e} de ses \'{e}l\'{e}ments donne lieu \`{a} au moins un \'{e}l\'{e}ment limite~\cite[p. 6]{fre2}.}
\end{dfn}

Since Fr\'{e}chet did not require that a compact set be closed, he defined the notion of an {\it extremal set}, which is closer to our modern day notion of compact.

\begin{dfn}
We shall call a set that is both closed and compact an {\it extremal} set; this nomenclature is justified later on. Within abstract set theory, extremal sets play a role akin to that of intervals in the theory of subsets of the real line. 
\end{dfn}


We might not know exactly why Fr\'{e}chet chose the word ``compact", but we have some evidence of why he did not choose some other terms. When Fr\'{e}chet first introduced the term, some mathematicians did not like his choice. For instance, Sch\"{o}nflies suggested that what Fr\'{e}chet called compact be called something like ``l\"{u}ckenlos'' (without gaps, closer to the modern notion of completeness) or ``abschliessbar'' (closable)~\cite[p. 266]{tay1}. Even to mature mathematicians the precise intuitions behind the term compactness was not yet clear.

Despite all of Fr\'{e}chet's early concern with intuitive definitions and
choice of terminology, it is surprising that at the end of his life, he could not
remember why he chose the term:
 
\begin{quote}
Doubtless I wanted to avoid a solid dense core with a single thread going off to
infinity being called compact. This is a hypothesis because I have completely
forgotten the reasons for my choice!~\cite[p. 440]{pier}.\footnote{From the original French: ... j'ai voulu sans doute \'{e}viter qu'on puisse appeler compact un noyau solide dense qui n'est agr\'{e}ment\'{e} que d'un fil allant jusqu'\`{a} l'infini. C'est une supposition car j'ai compl\`{e}tement oubli\'{e} les raisons de mon choix!}
\end{quote}

So even in the lifetime of the mathematician who named the concept, the original
intuition behind the concept was somewhat lost, and Fr\'{e}chet's intuitive nested
intersection definition was supplanted by less intuitive but more powerful notions of
limit-point, sequential, and open-cover compactness.

\subsection{Hausdorff: Compactness on metric spaces}
One of the obstacles to defining compactness, as we know it today, was to define it
in a way that would work for abstract topological spaces. This was a problem
for Fr\'{e}chet, and in the end he had to restrict his definition to V and E class, leaving open the question of defining compactness for L-classes, the abstract topological spaces. In 1914, Hausdorff made progress on this problem, introducing what we now call Hausdorff spaces, in which distinct points have distinct neighborhoods. In~\cite{haus} he introduced the term metric spaces (E spaces to Fr\'{e}chet), and he defined a set E to be compact if every infinite subset of E has a limit point in E, where limit point in this context means that every neighborhood of the point contains infinitely many
elements of E. Hausdorff's notion of compactness, which we would call limit-point compactness 
and is equivalent to countable compactness for Hausdorff spaces, remained the standard notion of compactness throughout the rapid development of point-set topology in the 1920s.\footnote{Thanks to an anonymous reviewer for this information and formulation.}

\subsection{Alexandroff and Urysohn: Open-cover compactness}

While Fr\'{e}chet was the first to formally define compactness, his contemporaries in
Russia, Alexandroff and Urysohn, appear to be the first to state it in its most general
form (in the context of abstract topological spaces). It is perhaps for this reason that the 
two Russians are often credited with defining the notion (e.g. \cite[p. 425S]{math}). Alexandroff and Urysohn were actually in close contact with Fr\'{e}chet~\cite[pp. 319-357]{tay2}.\footnote{Urysohn died, tragically, at the age of 26 in a swimming accident off the coast of France. Much of his work was published posthumously by Alexandroff, who kept up his correspondence with Fr\'{e}chet after Urysohn's untimely death.}

Also, though Alexandroff and Urysohn usually get credit for defining open-cover
compactness, Fr\'{e}chet was not unaware of the possibility of using neighborhoods to
characterize compactness. In a correspondence in 1905, Hadamard, Fr\'{e}chet's advisor, suggested that
he think in terms of neighborhoods to generalize the properties of to abstract topological spaces.\footnote{See \cite[p. 212]{koet}. The quotation comes from an undated letter from Hadamard to Fr\'{e}chet and can be found in full in \cite[p. 246]{tay1}).} The first definition that Fr\'{e}chet gave, in terms
of nested intersections, is the dual of (and hence logically equivalent to) countable
open-cover compactness.


\subsection{Open-cover vs. limit-point compactness}

Though Fr\'{e}chet may have been motivated originally to define compactness for
abstract topological spaces, he in fact restricted himself to metric spaces.  His approach
of looking at sequences and limits was not as general as the approach of using open covers,
which resulted in what we now take to be the correct definition of compactness.  Here we look 
at examples which illustrate why sequential compactness and open-cover compactness are not equivalent. We will use the concept of the least upper bound property, namely that any non-empty set containing an upper bound necessarily has a least upper bound.

\subsubsection{Sequentially compact does not imply compact}
Consider $S_\Omega =\{ \alpha \:|\: \alpha$ is an ordinal number and $\alpha < \Omega
\}$  with the order topology, where
$\Omega$ is the first uncountable ordinal number.  (See diagram below.  The first
{\em infinite} ordinal,
$\omega$, is the first ordinal after ``exhausting'' the natural numbers.  The first
{\em uncountable} ordinal,
$\Omega$, is the ordinal after ``exhausting" the countable ordinals.)\\

\begin{figure}[h]
\label{fig5}
 � \begin{center}
 � � \includegraphics[width=3in]{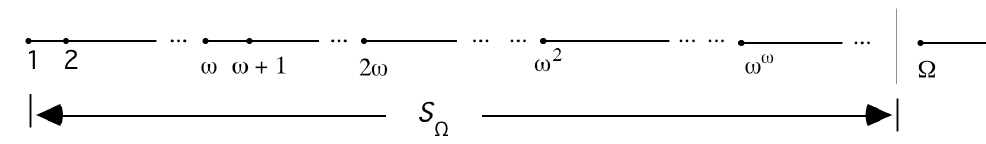}
 � \end{center}
 � \caption{Representation of $S_\Omega$.}
  \end{figure}
We know that all closed subsets of compact sets are compact (and all compact sets are closed).  So
$S_\Omega$  is not compact since it is not closed in the compact set $S_\Omega
\cup \{\Omega\}$.

However, it turns out that $S_\Omega$  is limit-point compact.  To see why this is
true, we will use the fact that any countable subset of $S_\Omega$ has an upper bound
in $S_\Omega$.  If we take any infinite subset of $S_\Omega$, it has a
countably infinite subset, which we will call $X$.  Since $X$ is countable, it has an
upper bound, let's call it $b$, in $S_\Omega$.  But the interval $[1,b]$  is compact since $S_\Omega$  has the l.u.b. property.  So there must be  a point in $[1,b]$  which is a limit point (of both $X$  and any set containing it).  Thus,  $S_\Omega$ is limit-point compact.  Essentially the same
argument shows that any sequence in $S_\Omega$  must have a convergent subsequence in
$S_\Omega$, so $S_\Omega$ is sequentially compact.

\subsubsection{Compact does not imply sequentially compact}
Just as we can have a space that is compact but not sequentially compact, we can also
have a space that is sequentially compact but not compact.\footnote{Thanks to an anonymous
reviewer for this example.} Consider the set of all functions from the interval $[0,1]$ to itself with 
the topology of pointwise convergence. This can be thought of as the infinite product $[0,1]^{[0,1]}$ 
with the product topology, which is compact by Tychonoff's theorem.\footnote{The product of compact
topological spaces is compact, or in German, Das Produkt von bikompakten R\"{a}umen ist wieder bikompakt, originally proved in \cite[p. 772]{tych}, though the theorem is sometimes credited to \v{C}ech \cite{foll}.} However, if $f_n(x)$ is the $n$th digit in the base-2 decimal expansion of $x$ (using the expansion that terminates in 0's if $x$ is a dyadic rational),
the sequence $f_n$, which is a sequence in the set of all functions from $[0,1]$ to itself, has no pointwise
convergent subsequence. It does have convergent subnets,\footnote{For example consider all functions which
map $[0,1]$ to the set $\{ 0,1 \}$. This will be a limit point for some subnets.} a concept that will be defined in the next section,
but not proper convergent sequences. 

\begin{figure}[h]
\label{fig6}
 � \begin{center}
 � � \includegraphics[width=3in]{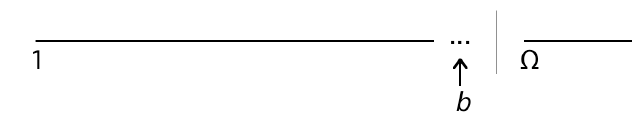}
 � \end{center}
 � \caption{Illustration of $b$ in $S_\Omega$.}
  \end{figure}


\section{Nets and Filters}

In the previous section we saw that the two important properties of compactness, those
stemming from the Bolzano-Weierstrass property (sequential compactness) and the Borel-Lebesgue property
(open-cover compactness), are not equivalent in abstract topological spaces. Open-cover compactness is more 
general and applicable, and for these reasons is considered to be compactness (and hence bears its name). However, it is possible to define open-cover compactness in a way that is analogous to sequential compactness, using modern notions of nets and filters, which we will develop here.\footnote{In this section we are more interested in the concepts and not the historical development, so we will be less careful than in earlier sections about giving original definitions and theorems, though we provide references for anyone who wants to track down the original formulations.} These two concepts are very different, on the surface, but they give rise to the same notion of convergence in abstract topological spaces.

\subsection{Moore and Smith: Nets}

The theory of nets was developed by E. H. Moore and his student H. L. Smith.\footnote{Little
biographical information about Smith is available. He received his Ph.D. from University of Chicago under Moore and got a job at Louisiana State University, but apparently after his important work on nets and filters, he dropped into obscurity in 1922 \cite[p. 563]{hal1}.} It is unclear whether Moore and Smith knew how nets could be used to define compactness. This connection is usually credited to Birkhoff  \cite[p. 64]{kell}, who applied Moore-Smith theory to general topological spaces, but Moore and Smith did generalize some of Fr\'{e}chet's compactness results in the same paper in which they defined nets \cite[p. 118]{moor}. Our goal here is to express compactness in terms of nets, so we will use the  
$S_\Omega$ example to motivate and illustrate net compactness.



The problem in the $S_\Omega$ example is that while $\Omega$ is a limit point of
$S_\Omega$  (any neighborhood of  $\Omega$ contains points of  $S_\Omega$), no
sequence in  $S_\Omega$ converges to $\Omega$ \cite[p. 76]{kell}. If we are limited to
taking a countable number of elements in the sequence, we will never reach $\Omega$. 
Nets provide one way of getting around this problem by allowing us to have something
like uncountable sequences. In our discussion of both nets and filters, we will consider only topological spaces, on which the notion of neighborhood is defined.\footnote{This treatment follows \cite[pp. 62-70]{kell} and
\cite[pp. 281-283, 286-289]{mcsh}.} 

To see how nets are a generalization of sequences, it is useful to think of
 sequences as functions on the natural numbers.

\begin{dfn}
A {\em sequence}, (denoted $\{x_n\}_{n \in \mathbb{N}} = \{x_1, x_2, x_3,\ldots\}$) is a function
which assigns to each element $n$ of the
 natural numbers, $\mathbb{N}$, a functional value $x_n$  in a set $X$.  
\end{dfn}

We would like to replace $\mathbb{N}$ with a set that can be uncountable but has an ordering
similar to that of $\mathbb{N}$. In other words, we want to stipulate conditions for an
ordering relation on a generic set that generalizes the way $>$ orders natural
numbers. We will call this relation $\succ$  to suggest the connection to $>$, and we
will say that this relation ``directs" a given set.

\begin{dfn}
A non-empty set $D$, with the relation $\succ$ is called {\em directed} iff
\begin{tabular}{r@{  }l}
 (i)  &if $d_1, d_2, d_3 \in D$  such that $d_1 \succ d_2$  and
$d_2 \succ d_3$  then $d_1 \succ d_3$;\\
 (ii)  &if  $d_1, d_2 \in D$, then there is a $d_3 \in D$  such that $d_3
\succ d_1$  and $d_3 \succ d_2$.\\
\end{tabular}
\end{dfn}

So the definition of a net is simply the definition of a sequence with $\mathbb{N}$  replaced
 by the notion of a directed set. From now on, $D$  will stand for a directed set with
the relation $\succ$  as defined above. 

\begin{dfn}
A net (denoted $\{x_d\}_{d \in D}$  or simply $\{x_d\}$)  is a function which assigns
to each element $d$  of
 a directed set $D$  a functional value $x_d$  in a set $X$.
\end{dfn}

Once we know what a net is, we can state what it means for it to converge. Again we
can derive the definition for net convergence and limit point by taking the
definitions involving sequences and simply replacing $\mathbb{N}$   and $>$ with $D$ and 
$\succ$.

\begin{dfn}
A net $\{x_d\}$ {\em converges} to $a \in X$  (denoted $\{x_d\} \rightarrow a$) iff for
every neighborhood $U$  of  $a$,
 there is an index $d_0 \in D$  such that if $d \succ d_0$   then $x_d \in U$ (i.e. if
the net is {\em eventually in} each neighborhood of  $a$). 
\end{dfn}

\begin{dfn}
A point $a$ is a {\em limit point} of $\{x_d\}$ if for every neighborhood $U$ of $a$ 
and every $d_0 \in D$ there is a $d \succ d_0$ such that $x_d \in U.$
\end{dfn}

In order to state compactness in terms of nets, we also need the concept of subnet,
the analog of subsequence.\footnote{Incidentally, Kelley, who first coined the term 
``net" had considered using the term ``way" so the analog of subsequence would be
``subway."  McShane also proposed the term ``stream" for net since he thought it was
intuitive to think of the relation of the directed set as ``being downstream from" \cite[p. 282]{mcsh}.}   
Part of the definition of subsequence generalizes easily, but the other part requires 
us to think about subsequences in a slightly different way than we are accustomed. 
The first defining property of subsequence is that each
element of the subsequence can be identified with an element of the sequence. This
property is generalized in (i) below. The second defining property requires that the
subsequence is ordered in a similar way as the sequence. Usually we require the
indices of the subsequence, like the indices of the sequence, to be strictly
increasing. In other words, for a subsequence $\{x_{n_k}\}$  of a sequence $\{x_n\}$,
the $n_k$  are positive integers such that $n_1< n_2<n_3 \cdots$. But the feature of
this condition which turns out to be important is simply the fact that as $k
\rightarrow \infty$, so do the $n_k$. This property is generalized in (ii) below.

\begin{dfn}
A {\em subnet} of a net $\{x_d\}_{d \in D}$  is a net $\{y_b\}_{b \in B}$  where $B$
is a directed set and there is a function $\varphi:B \rightarrow D$  such that:\\
\begin{tabular}{r@{  }l}
 (i) & $y_b = x_{\varphi(b)}$      and \\
 (ii) & $\forall  d \in D, \exists b_0 \in B$ such that if $b \succ b_0$ then
$\varphi(b) \succ d$. \\
\end{tabular}
\end{dfn}

We are now ready to characterize compactness in terms of nets.

\begin{thm}
A topological space $X$  is {\em compact} iff either \ldots
\newline $\bullet$  Every net of points of $X$  has a limit point in  $X$, or
\newline $\bullet$ Every net of points of $X$  has a convergent subnet in $X$.

\end{thm}

Notice that these definitions are precisely the same as limit point and sequential
 compactness for metric spaces with the term ``net" substituted for ``sequence."

Applying these definitions to the $S_\Omega$  example, we can show why $S_\Omega$ is
not compact.
  If we take a net ${\{x_d\}}$ of elements of  $S_\Omega$, it is no longer the case that
there will necessarily be a limit in $S_\Omega$. In particular, let $D = S_\Omega$ 
and 
$x_d = d$. Then ${\{x_d\}}$  converges to $\Omega$, which is not in $S_\Omega$. Thus no
subnet of ${\{x_d\}}$ will converge to a point in $S_\Omega$. 

\subsection{Cartan (and Smith): Filters}

Nets are not the only way of generalizing sequences. Another generalization of
sequence is a filter, a notion suggested by Cartan in 1937.\footnote{See also \cite[p. 8]{bou2}.} While different from a net, both nets
and filters give rise to the same notion of convergence on topological spaces. That is to say, on 
abstract topological spaces, they are essentially the same.\footnote{To show equivalence on abstract topological spaces, there is for example an exercise in Kelley that establishes a 
dictionary mapping between them (i.e. given a net you can find a filter, and vice
versa)\cite[p. 83]{kell}. But there is a subtle distinction for a particular type of
limit found in the advanced theory of integration \cite[p. 371]{smit}.} Nonetheless, some mathematicians find nets more intuitively appealing and useful, while others prefer filters.

The idea behind filters was foreshadowed by Riesz \cite{ries} in 1907 when he provided axioms for topology based on limit
points instead of metrics, though his topological axioms are not equivalent to the standard
ones we use today, and his work did not result in a fruitful line of research. Riesz defines a concept 
called an ``ideal"  which is essentially the same as what we now call an ultrafilter. Smith independently discovered filters as an attempt to explain what was lacking in the theory of nets
that he and Moore proposed. 

Again, we will define the notions we need to state compactness in terms of filters and
then apply our compactness result to show  $S_\Omega$ is not compact. As with nets,
we can look at convergence of sequences to motivate the idea of convergence of
filters. However, whereas with nets the focus was on the index set, with filters the
focus is on neighborhoods.\footnote{This treatment follows \cite{dixm}.}  

\begin{dfn}
Let $X$ be a set. A set $\Phi$  of subsets of $X$ is called a {\em filter} iff
\begin{tabular}{r@{  }l}
  (i)  & $\emptyset \notin \Phi$\\
  (ii) &  $A_1 \subset A_2 \subset X$ and $A_1 \in \Phi \Rightarrow A_2 \in
\Phi$\\
  (iii) &  $A_1, A_2 \in \Phi \Rightarrow A_1 \cap A_2 \in \Phi$
\end{tabular}\\
\end{dfn}

As with nets, we should define what it means for a filter to converge:

\begin{dfn}
A filter $\Phi$  converges to $a \in A$ (denoted $\Phi \rightarrow a$)  iff each
neighborhood of $a$  is a member of  $\Phi$.

\end{dfn}

There is a natural way to associate a filter with any sequence. If $x_1, x_2, x_3,
\ldots$  is a sequence in $X$, we can associate with this sequence a filter $\Phi$ on
$X$ such that  $\forall a \in X, \{x_n\} \rightarrow a$ iff $\Phi \rightarrow a$. In
particular, let $\Phi=\{A \subset X \:|\: \exists k_A$ such that $\forall i\geq k_A,
x_i \in A\}$. So the tails of the sequence are contained in neighborhoods which are members of the
filter. The condition that each neighborhood of $a$ is in the filter is then
equivalent to the condition that the sequence is eventually in any neighborhood of
$a$.

\begin{figure}[h]
\label{fig7}
 � \begin{center}
 � � \includegraphics[width=3.5in]{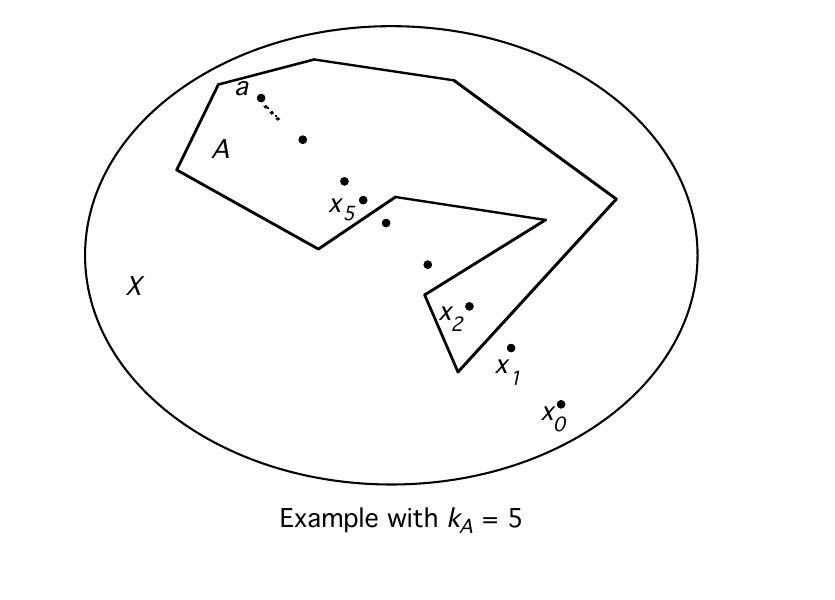}
 � \end{center}
 � \caption{Member of filter with $k_A = 5$.}
  \end{figure}

In order to define compactness in terms of filters, we need one more notion, that of
an ultrafilter.

\begin{dfn}
A filter in $X$ is an {\em ultrafilter} iff no filter in $X$ properly contains it. 
\end{dfn}

The notion of ultrafilter is not exactly analogous to subsequence, but in the
formulation of compactness, it serves the same purpose.

\begin{thm}
 A topological space is {\em compact} iff every ultrafilter on  $X$ converges
in $X$.
\end{thm}

Now we can return again to our example and get a sense in terms of filters for why  
$S_\Omega$ is not compact.\footnote{The proof of this claim, in particular when we
assert that there is an ultrafilter containing our filter, actually relies on the
axiom of choice.}  We want to show that there is an ultrafilter on $S_\Omega$  that
does not converge. Consider all the neighborhoods of $\Omega$  in  $S_\Omega \cup
\Omega$. Let $\Phi=\{A \subset S_\Omega \:|\: \exists \alpha \in S_\Omega$ such that
$\forall \beta \geq \alpha, \beta \in A\}$.  

\begin{figure}[h]
\label{fig8}
 � \begin{center}
 � � \includegraphics[width=3in]{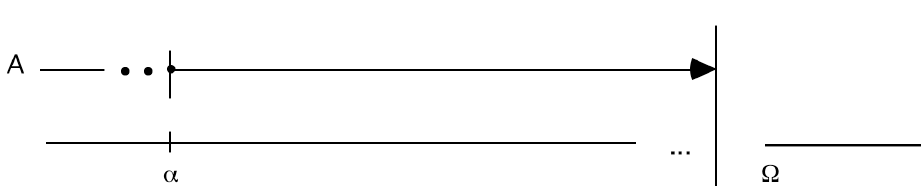}
 � \end{center}
 � \caption{Member of filter on $S_\Omega$.}
  \end{figure}
This clearly satisfies the definition of a filter. Let $\Psi$ be any ultrafilter
containing $\Phi$. We claim $\Psi$ does not converge in $S_\Omega$. Suppose it did. 
Say that $\Psi \rightarrow b$. Now pick some $\alpha > b, \alpha \in S_\Omega$. 
Then $A^+ = \{\beta:\beta \geq \alpha\} \in \Psi$  since $A^+ \in \Phi \subseteq
\Psi$. We also have $A^- = \{\beta:\beta < \alpha\} \in \Psi$ since $A^-$  is an open
neighborhood of $b$ (and we claim $\Psi$ converges to $b$).  

\begin{figure}[h]
\label{fig9}
 � \begin{center}
 � � \includegraphics[width=3in]{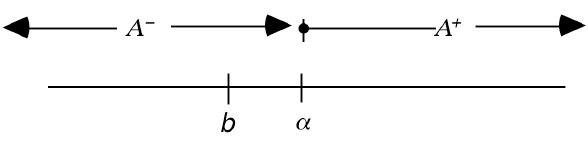}
 � \end{center}
 � \caption{Illustration of $A^-$ and $A^+$.}
  \end{figure}
But $A^+ \cap A^- = \emptyset$, which violates the definition of a filter, so our
assumption must be wrong. Thus, $\Psi$ must not converge, and hence $S_\Omega$  is not
compact.


\subsection{Final Comments}
Here ends the story, a sort of co-evolution of the two different, but related notions of sequential and open-cover compactness. Today when we use the term ``compact'' we mean open-cover compact, but this paper, and the terms listed in the Appendix shows that this was not always the case. The story of how open-cover compactness came to be seen as the {\em right} one is a story of developing mathematics without always knowing where it going, how important terms should be defined, and how widely they might be applied. 

It might be worth noting, in closing, that even this paper, which attempts to characterize the evolution of compactness is only a sort of snapshot. In the time it took to write this paper, and revise it, textbooks have changed. For instance, the latest version of a standard general topology textbook~\cite{mun2} now includes a discussion of nets and filters, whereas the earlier version~\cite{mun1}, available during the writing of the original version of this paper, did not.

The lesson to be drawn is simply that mathematics evolves and changes as concepts become clearer and are applied in more general situations. This may be obvious to the mathematician who is involved in making these conceptual advances, but may be less clear to the student who sees a textbook as a collection of established facts. Textbooks are, as perhaps they should be, a distillation of what we currently know. They are also historical documents, in their own right, and being aware of that fact may help students mature as mathematicians.




\pagebreak
\pagebreak
\subsection{Appendix :  Terminology}
There are many notions related to (but not necessarily equivalent to) compactness. Table~\ref{t:flavors} contains a list of some of these notions.\\


\begin{table}[h]
\centering

\begin{tabular}{|p{4.3 in}|}

\hline
[open-cover] \emph{compact}:  Every open cover has a finite subcover.\\
(also called the Borel-Lebesgue property)\\\

\emph{sequentially compact}:  Every sequence has a convergent subsequence.\\\

\emph{countably compact}:  Every countable open cover has a finite subcover.\\\

\emph{limit-point compact}:  Every infinite subset of $X$  has a limit point in $X$.\\
(also called Fr\'{e}chet compact or the Bolzano-Weierstrass property)\\\

\emph{relatively compact}:  The closure is compact.\\\

\emph{quasi-compact}:  Compact to Bourbaki.\\\

\emph{pseudo-compact}:  Each continuous real valued function on $X$  is bounded.\\\

\emph{finally compact}:  index of compactness is ${\aleph_0}$.\\
(also called Lindel\"{o}f compact)\\

\hline

\end{tabular}

\caption{Flavors of compactness}
\label{t:flavors}

\end{table}

Many of these concepts are related. For instance, compactness
implies countable compactness implies limit point compactness. Sequential compactness
implies countable compactness.  And if we put further restrictions on our spaces we
can get implications in the other direction. In $T_1$ spaces, limit-point compactness
implies countable compactness. In first countable spaces, countable compactness
implies sequential compactness. In second countable spaces, sequential compactness
implies compactness. In particular, we know that in compact metric spaces, which turn out to
be second countable, the first four notions of compactness in Table~\ref{t:flavors}
are equivalent.

It took some time as compactness was applied to different types of spaces for
relationships like these to be worked out. It also took time for names to stabilize. Table~\ref{t:compactnames} lists different terms used for compactness-related ideas used by some of the most influential mathematicians in the historical development.
In this paper, the modern names have been used.

\pagebreak
\begin{table}[h]
\centering
\begin{tabular}{|p{1.1 in}|p{0.45 in}|p{0.9 in}|p{1.35 in}|}
\hline

\bf{Who} & \bf{When}  & \bf{Their term}  & \bf{Modern
term}\\

\hline Fr\'{e}chet   & 1906       & compact    & relatively sequentially compact\\
\cline{3-4}
              &            & extremal         & sequentially compact\\ \hline
Russian School \newline (Alexandroff, etc.)  & 1920's  & bicompact  & compact\\
\cline{3-4}
              &            & compact         & countably \newline compact\\
\hline          Bourbaki      & 1930's     & quasi-compact  & compact\\ \cline{3-4}
              &            & compact         & compact and \newline Hausdorff\\
\hline      

\end{tabular}
\caption{Names of historical compactness-related terms}
\label{t:compactnames}
\end{table}


\end{document}